\documentclass[12pt]{article}
\usepackage{amsmath}
\usepackage{amssymb}
\usepackage{eucal}
\usepackage{amsthm}
\usepackage{graphicx,graphics}
\numberwithin{equation}{section}
\usepackage{pdfsync}
\def \dis {\displaystyle}

\def \confai {\;-\kern -.5em\rightharpoonup\;}
\def \MMMMI#1{-\mskip -#1mu\int}
\def \moy {\displaystyle\MMMMI{19.2}}

\def \cqfd {\hfill\Box}
\def \al {\alpha}

\def \Ga {\Gamma}

\def \De {\Delta}
\def \ep {\varepsilon}

\def \ph {\varphi}
\def \th {\theta}
\def \ka {\kappa}
\def \si {\sigma}
\def \Si {\Sigma}

\def \ZZ {\mathbb Z}

\def \RR {\mathbb R}

\def \Y {\mathcal{Y}}
\def \beq {\begin{equation}}
\def \eeq {\end{equation}}
\def \ba {\begin{array}}
\def \ea {\end{array}}
\def \bs {\bigskip}
\def \ms {\medskip}

\def \ecart {\noalign{\medskip}}

\newtheorem{Thm}{Theorem}[section]

\newtheorem{Adef}[Thm]{Definition}

\newtheorem{Arem}[Thm]{Remark}
\newenvironment{Rem}{\begin{Arem}\rm}{\end{Arem}}
\newtheorem{Aexa}[Thm]{Example}

\newtheorem{Anot}[Thm]{Notation}

\def \refe #1.{(\ref{#1})}
\def \reff #1.{figure~\ref{#1}}
\def \refs #1.{Section~\ref{#1}}
\def \refss #1.{subsection~\ref{#1}}
\def \refD #1.{Definition~\ref{#1}}
\def \refT #1.{Theorem~\ref{#1}}
\def \refL #1.{Lemma~\ref{#1}}
\def \refC #1.{Corollary~\ref{#1}}
\def \refP #1.{Proposition~\ref{#1}}
\def \refPt #1.{Properties~\ref{#1}}
\def \refR #1.{Remark~\ref{#1}}
\def \refE #1.{Example~\ref{#1}}
\def \refN #1.{Notation~\ref{#1}}
\title{Bounds on strong field magneto-transport in three-dimensional composites}
\author{\begin{tabular}{cc}
    {Marc BRIANE }                                             & {Graeme W. MILTON}                            			\\*[-.0em]
    {\small  Institut de Recherche Math\'ematique de Rennes} & {\small Department of Mathematics}	\\*[-.3em]
    {\small  Universit\'e Europ\'eenne de Bretagne} & {\small University of Utah} 									\\*[-.3em]
    {\small mbriane@insa-rennes.fr}                            & {\small milton@math.utah.edu}
\end{tabular}}
\begin{document}
\maketitle
\begin{abstract}
This paper deals with bounds satisfied by the effective non-symmetric conductivity of three-dimensional composites in the presence of a strong magnetic field. On the one hand, it is shown that for general composites the antisymmetric part of the effective conductivity cannot be bounded solely in terms of the antisymmetric part of the local conductivity, contrary to the columnar case studied in \cite{BM3}. So, a suitable rank-two laminate the conductivity of which has a bounded antisymmetric part together with a high-contrast symmetric part, may generate an arbitrarily large antisymmetric part of the effective conductivity. On the other hand, bounds are provided which show that the antisymmetric part of the effective conductivity must go to zero if the upper bound on the antisymmetric part of the local conductivity goes to zero, and the symmetric part of the local conductivity remains bounded below and above. Elementary bounds on the effective moduli are derived assuming
the local conductivity and effective conductivity have transverse isotropy in the plane orthogonal to the magnetic field.  New Hashin-Shtrikman type bounds for two-phase three-dimensional composites with a non-symmetric conductivity are 
provided under geometric isotropy of the microstructure. The derivation of the bounds is based on a particular variational principle symmetrizing the problem, and the use of $Y$-tensors involving the averages of the fields in each phase.
\end{abstract}
\section{Introduction}
It is known since the seminal discovery of Hall in the end of the 19th century \cite{Hall}, that a low magnetic field perturbs the matrix resistivity (or equivalently the conductivity) of a conductor by inducing a small non-symmetric part characterized by the so-called Hall coefficient. In the 80's Bergman \cite{Ber} first gave a general formula for the effective Hall coefficient involving currents that solve the symmetric conductivity equations in the absence of a magnetic field. However, there are few explicit formulas for the effective Hall coefficient except in very particular cases like two-phase two-dimensional composites \cite{MilHall,BerLiSt,BMM}, or columnar composites \cite{BerSt1,BerSt2,StrelBer,Grab1,Grab2}. The situation is still less favorable in the strong field case \cite{BerSt1a,BerSS}, namely when the symmetric part and the antisymmetric part of the conductivity are of the same order. In three dimensions, only when the antisymmetric
part is constant do we have an exact formula for the antisymmetric part of the effective tensor
\cite{StBer}. So, rather than trying to get explicit relations for the effective tensors it seems more practical to derive bounds. The theory of the bounds in homogenization has gone through a considerable development since the original variational approach of Hashin and Shtrikman \cite{HaSh}. We refer to \cite{MilBook} for a comprehensive survey. In fact, very little is known about bounds for strong field magneto-transport. Recently, we derived in \cite{BM3} optimal bounds for multiphase columnar composites. The aim of this paper is to extend, at least partially, the result of \cite{BM3} to two-phase three-dimensional composites.
\par
In the present context we consider a three-dimensional conductor having a periodic structure (this is actually not a restrictive assumption) in the presence of a fixed vertical strong magnetic field. Under the transverse isotropy assumption the local conductivity of the conductor takes the general form
\beq\label{isocon}
\si(y)=\begin{pmatrix} a(y) & -c(y) & 0 \\ c(y) & a(y) & 0 \\ 0 & 0 & b(y) \end{pmatrix},\quad\mbox{for }y\in\RR^3,
\eeq
where the coefficient $c(y)$ is induced by the presence of the magnetic field parallel to the $y_3$-axis, which also influences $a(y)$ and $b(y)$ and causes them to be non-equal in the case of a conductor that is isotropic in the absence of the magnetic field.
Similarly, assuming an transversely isotropic microstructure, or at least one that is invariant under $90^\circ$ or $120^\circ$ rotations about the $y_3$-axis, the constant effective conductivity of the composite is given by
\beq\label{isocon*}
\si_*=\begin{pmatrix} a_* & -c_* & 0 \\ c_* & a_* & 0 \\ 0 & 0 & b_* \end{pmatrix}.
\eeq
Our goal is to derive bounds for the effective coefficients $a_*$, $b_*$, $c_*$ of $\si_*$ in terms of the coefficients $a(y)$, $b(y)$, $c(y)$ of the local  conductivity $\si(y)$.
\par
In Section \ref{s.elem}, we derive elementary bounds (see Theorem~\ref{thm.elem}) on the effective coefficients $a_*$, $b_*$ and $c_*$. These are obtained by taking uniform trial fields in 
a variational principle for non-symmetric tensors
deduced in \cite{Mil90,FaPa} from a symmetrization of the constitutive law $j=\si e$, and its adjoint, adapted from the variational approach performed in \cite{ChGi} for complex tensors.

\par
In Section \ref{s.cex}, we show (see Theorem~\ref{thm.bou}) that contrary to the columnar case of \cite{BM3}, it is not possible to bound the antisymmetric part of the effective conductivity $\si_*$ only in terms of the coefficients $c(y)$. Indeed, when $\si(y)$ is independent of $y_3$ for a vertical columnar structure, the key ingredient for the derivation of the optimal bounds in \cite{BM3} is based on the positivity of the $(2\times 2)$ determinant of the local electric field \cite{AlNe1,AlNe2}, i.e.
\beq\label{posdet}
\De_{1,2}(DU):=\partial_1 u_1\,\partial_2 u_2-\partial_1 u_2\,\partial_2 u_1>0\quad\mbox{a.e. in }\RR^3,
\eeq
where the vector-valued potential $U=(u_1,u_2,u_3)$ solves the conductivity problem
\beq\label{conpro}
\left\{\ba{ll}
{\rm div}\left(\si DU\right)=0 & \mbox{in }\RR^3
\\
U(y)-y & \mbox{is $Y$-periodic}.
\ea\right.
\eeq
Due to a suitably constructed rank-two laminate with high-contrast conductivity, we prove simultaneously that the inequality \refe{posdet}. does not hold, and that arbitrarily large effective coefficients $c_*$ can be obtained while the local coefficient $c(y)$ is bounded. This negative result agrees with the pathologies obtained in \cite{BM1,BM2} with different microstructures, related to bounds on the effective Hall coefficient in the low magnetic field regime.
As a consequence, a bound for $c_*$ involves both the upper bound for $|c(y)|$ and the bounds from below and above for $a(y)$ in \refe{isocon}. (see Theorem~\ref{thm.bou} and Remark~\ref{rem.isot}). This bound shows $c_*\to 0$
when the upper bound on $|c(y)|$ goes to zero, provided $a(y)$ remains bounded from below and above.
\par
In Section \ref{s.HSbou}, to improve the previous bounds we restrict ourselves to a two-phase local conductivity
\beq\label{isocon2}
\si(y)=\chi_1(y)\begin{pmatrix} a_1 & -c_1 & 0 \\ c_1 & a_1 & 0 \\ 0 & 0 & b_1 \end{pmatrix}
+\chi_2(y)\begin{pmatrix} a_2 & -c_2 & 0 \\ c_2 & a_2 & 0 \\ 0 & 0 & b_2 \end{pmatrix},
\quad\mbox{for }y\in\RR^3,
\eeq
with prescribed volume fraction $f_i=\langle\chi_i\rangle$, for $i=1,2$, with $f_1+f_2=1$. We then derive (see Theorem~\ref{thm.HSbou}) Hashin-Shtrikman type bounds for the effective conductivity $\si_*$, involving three intermediate coefficients $a_Y$, $b_Y$, $c_Y$ which are explicitly expressed in terms of the entries of $\si_i$, for $i=1,2,*$. In particular, it is shown that the point $(a_Y,-c_Y)$ belongs to a disk which is tangent to the axis $a=0$ at some point $(0,c)$, and which contains the disk passing through the points $(a_1,c_1)$, $(a_2,c_2)$, and tangent to the axis $a=0$ at the same point $(0,c)$ (see Figure~\ref{fig1} below). The derivation of these new bounds is based on a combination of three main ingredients:
\begin{itemize}
\item the geometric isotropy of the phases defined in \cite{Wil} for random composites,
\item the variational principle for non-symmetric tensors,
\item the use of $Y$-tensors similar to \cite{GiMi} (see also \cite{Berry}), giving relations between the averages of the fields in each phase.
\end{itemize}
\subsection*{Notations}
\begin{itemize}
\item $(e_1,e_2,e_3)$ denotes the canonic basis of $\RR^3$.
\item $I$ denotes the unit matrix of $ \RR^{3\times 3}$, and
$J:=\left(\begin{smallmatrix} 0 & -1 & 0 \\ 1 & 0 & 0 \\ 0 & 0 & 0 \end{smallmatrix}\right)$.
\item For any matrix $M\in\RR^{2\times 2} $, $M^T$ denotes the transpose of $M$, $M^S:=\frac{1}{2}\,\bigl(M+M^T\bigr)$ the symmetric part of $M$, and $M^A:=\frac{1}{2}\,\bigl(M-M^T\bigr)$ the antisymmetric part of $M$.
\item $Y$ denotes the unit cube $[0,1]^3$, and $\langle\cdot\rangle$ the $Y$-average.
\item For a function $f$ defined on the unit sphere $S^2$ of $\RR^3$, $\langle f\rangle_{S^2}$ denotes the average of $f$ over $S^2$, i.e.
\beq\label{moyf(xi)}
\langle f\rangle_{S^2}:=\moy_{S^2}f(\xi)\,d\xi={1\over 4\pi}\int_0^{2\pi}d\ph\int_0^\pi f(\sin\th \cos\ph,\sin\th \sin\ph,\cos\th)\,\sin\th\,d\th.
\eeq
\item For $\alpha,\beta>0$, ${\cal M}_\sharp(\alpha,\beta;Y)$ denotes the set of the $Y$-periodic invertible matrix-valued functions $A:\RR^3\to\RR^{2\times 2}$ such 
that\begin{gather}
\label{Mab}
\forall\,\xi\in {\RR^3,}\quad A(y)\,\xi\cdot\xi\geq \alpha\,|\xi|^2\quad\text{and}\quad A^{-1}(y)\,\xi\cdot\xi\geq \beta^{-1}\,|\xi|^2 \quad\text{a.e.~}y\in Y.
\end{gather}
\item $L^2_\sharp(Y)$ denotes the space of the $Y$-periodic functions, which are square integrable in $Y$.
\item $H^1_\sharp(Y)$ denotes the space of the $Y$-periodic functions, with gradient in $L^2_\sharp(Y)^3$.
\item For $u:\RR^3\longrightarrow\RR$, $\nabla u:=\bigl(\frac{\partial u}{\partial x_i}\bigr)_{1\leq i\leq 3}$.
\item For $U:\RR^3\longrightarrow\RR^3$, $U=(u_1,u_2,u_3)$, $DU:=\bigl(\frac{\partial u_j}{\partial x_i}\bigr)_{1\leq i,j\leq 3} $.
\item For $\Si:\RR^3\longrightarrow\RR^{3\times 3}$, $\mathrm{div}\left(\Si\right):=\bigl(\frac{\partial\Si_{ij}}{\partial x_i}\bigr)_{1\leq j\leq 3}$.
\end{itemize}
\section{Elementary bounds on magneto-transport}\label{s.elem}
To derive bounds on composites we may assume  that the associated microstructures are $Y$-periodic (see, e.g., \cite{All} Theorem~1.3.23), where $Y$ is any cube of $\RR^3$, say $Y=[0,1]^3$. 
In this section and the next we consider a three-dimensional $Y$-periodic conductor in the presence of a strong magnetic field parallel to the $y_3$-axis so that the resulting matrix-valued conductivity $\si(y)$ is given by
\beq\label{iso.sig}
\si(y)=\begin{pmatrix}
a(y) & -c(y) & 0
\\
c(y) & a(y) & 0
\\
0 & 0 & b(y)
\end{pmatrix},
\quad\mbox{for }y\in\RR^3,
\eeq
where the coefficients $a(y),b(y),c(y)$ satisfy for prescribed positive numbers $\underline{a},\overline{a},\overline{c}>0$, with
$\underline{a}\leq\overline{a}$, the following bounds
\beq\label{bou.abc}
\underline{a}\leq a(y)\leq\overline{a}\quad\mbox{and}\quad|c(y)|\leq\overline{c},\qquad\mbox{a.e. }y\in \RR^3.
\eeq
By virtue of the periodic homogenization formula (see, e.g., \cite{BLP}) the effective conductivity $\si_*$ associated with $\si(y)$ is given by
\beq\label{sig*}
\si_*=\langle\si DU\rangle,\quad\mbox{where the potential $U$ solves }\quad
\left\{\ba{ll}
{\rm div}\left(\si DU\right)=0 & \mbox{in }\RR^3
\\
U(y)-y & \mbox{is $Y$-periodic}.
\ea\right.
\eeq
Recall that $\si_*$ is also the homogenized conductivity obtained from the oscillating sequence $\si({x\over\ep})$ as $\ep\to 0$ by a homogenization process (see, e.g., \cite{BLP}).

Now consider a periodic electric field $e\in L^2_\sharp(Y)^3$ and a periodic current field $j\in L^2_\sharp(Y)^3$ that solves the conductivity equations
\beq\label{con1}
j=\si e,\quad  {\rm div}~j=0,\quad {\rm curl}~e=0
\eeq
and another periodic electric field $e'\in L^2_\sharp(Y)^3$ and another periodic current field $j'\in L^2_\sharp(Y)^3$ that solves the adjoint equations
\beq\label{con2}
j'=\si^T e',\quad  {\rm div}~j'=0,\quad {\rm curl}~e'=0.
\eeq
The average fields are related by the effective tensor $\si_*$:
\beq\label{efften} 
\langle j\rangle=\si_*\langle e\rangle,\quad \langle j'\rangle=\si_*^T\langle e'\rangle.
\eeq
Define the symmetric tensor
\beq\label{L(y)}
L(y):=\begin{pmatrix}
(\si^S)^{-1} & -\,(\si^S)^{-1}\,\si^A
\\ \ecart
\si^A\,(\si^S)^{-1} & \si^S-\si^A\,(\si^S)^{-1}\,\si^A
\end{pmatrix}(y).
\eeq
Then, an easy computation yields
\beq\label{LEF}
F=\begin{pmatrix} e_S \\ j_A \end{pmatrix}=L\begin{pmatrix} j_S \\ e_A \end{pmatrix}=LE,\quad\mbox{where}\quad
\left\{\ba{ll}
\dis e_S:={1\over 2}\left(e+e'\right), & \dis e_A:={1\over 2}\left(e-e'\right)
\\ \ecart
\dis j_S:={1\over 2}\left(j+j'\right), & \dis j_S:={1\over 2}\left(j+j'\right).
\ea\right.
\eeq
Moreover, mimicking the approach of \cite{ChGi} for complex tensors, extended in \cite{Mil90,FaPa} (see also \cite{MilBook}, p.~277) for real but non-symmetric tensors, the following variational principle holds
\beq\label{varpri}
\begin{pmatrix} j_0 \\ e_0 \end{pmatrix}^{\!\! T}\!\! L_*\begin{pmatrix} j_0 \\ e_0 \end{pmatrix}=
\min\left\{
\left\langle\! \begin{pmatrix} j_S \\ e_A \end{pmatrix}^{\!\! T}\!\! L\begin{pmatrix} j_S \\ e_A \end{pmatrix}\!\right\rangle:
{\left|\,\ba{lll} e_A\in L^2_\sharp(Y)^3, & \!\! {\rm curl}\left(e_A\right)=0, & \!\!\langle e_A\rangle=e_0
\\ \ecart
j_S\in L^2_\sharp(Y)^3, & \!\! {\rm div}\left(j_S\right)=0, & \!\! \langle j_S\rangle=j_0.\ea\right.}\!\!\right\},
\eeq
with the symmetric effective tensor
\beq\label{L*}
L_*:=\begin{pmatrix}
(\si_*^S)^{-1} & -\,(\si_*^S)^{-1}\,\si_*^A
\\ \ecart
\si_*^A\,(\si_*^S)^{-1} & \si_*^S-\si_*^A\,(\si_*^S)^{-1}\,\si_*^A
\end{pmatrix}.
\eeq
By substituting constant trial fields $e_A=e_0$ and $j_S=j_0$ in the variational principle one immediately obtains the elementary bound
\beq\label{elem}
L_*\leq \langle L\rangle.
\eeq
This elementary bound implies the following theorem:
\begin{Thm}\label{thm.elem}
Assuming $\si_*$ and $\si(y)$ have the forms \refe{isocon*}. and \refe{iso.sig}. the constant $b_*$ must satisfy the arithmetic and harmonic mean bounds
\beq\label{elemb}
1/\langle 1/b \rangle\leq b_*\leq\langle b \rangle,
\eeq
and the pair $(a_*,c_*)$ must satisfy the circle bounds (which confine $(a_*,c_*)$ to lie within a circle in the $a_*$-$c_*$ plane)
given by
\beq\label{simpcirc} 
(c_*-c_L)^2\leq\left(a_*-a_L\right)\left(d_L-a_*\right),
\eeq
where 
\beq\label{siL}
a_L:=\left\langle{1\over a}\right\rangle^{-1},\quad c_L:=\left\langle{c\over a}\right\rangle a_L,\quad
d_L:=\left\langle a+{c^2\over a}\right\rangle-{c_L^2\over a_L}.
\eeq
\end{Thm}
\noindent
\begin{Rem}\label{rem.alt}
Taking the minimum in \refe{varpri}. over all fields $e_A$ and $j_S$ with $\langle e_A\rangle=e_0$ and $\langle j_S\rangle=j_0$, and ignoring the differential constraints that
${\rm curl}\left(e_A\right)=0$ and ${\rm div}\left(j_S\right)=0$ gives the elementary bound $L_*^{-1}\leq \langle L^{-1}\rangle$. However this does not yield any new inequalities
beyond \refe{elemb}. and \refe{simpcirc}. due to the structure of the matrices $L_*$ and $L(y)$.
\end{Rem}
\noindent
{\bf Proof of Theorem~\ref{thm.elem}.}
The proof follows the proof of the elementary bounds in Proposition~3.1 of \cite{BM3}.
Assuming $\si_*$ and $\si(y)$ have the forms \refe{isocon*}. and \refe{iso.sig}. we obtain
\beq\label{L*big}
L_*=\begin{pmatrix}
\dis {1\over a_*} & 0 & 0 & 0 & \dis {c_*\over a_*} & 0
\\
0 & \dis {1\over a_*} & 0 & \dis -\,{c_*\over a_*} & 0 & 0
\\
0 & 0 & \dis {1\over b_*} & 0 & 0 & 0
\\
0 & \dis -\,{c_*\over a_*} & 0 & \dis a_*+{c_*^2\over a_*} & 0 & 0
\\
\dis {c_*\over a_*} & 0 & 0 & 0 & \dis a_*+{c_*^2\over a_*} & 0
\\
0 & 0 & 0 & 0 & 0 & b_*
\end{pmatrix}.
\eeq
and
\beq\label{aveL}
\langle L\rangle
=\begin{pmatrix}
\dis {1\over a_L} & 0 & 0 & 0 & \dis {c_L\over a_L} & 0
\\
0 & \dis {1\over a_L} & 0 & \dis -\,{c_L\over a_L} & 0 & 0
\\
0 & 0 & \left\langle\dis {1\over b}\right\rangle & 0 & 0 & 0
\\
0 & \dis -\,{c_L\over a_L} & 0 & \dis d_L+{c_L^2\over a_L} & 0 & 0
\\
\dis {c_L\over a_L} & 0 & 0 & 0 & \dis d_L+{c_L^2\over a_L} & 0
\\
0 & 0 & 0 & 0 & 0 & \langle b \rangle
\end{pmatrix}.
\eeq
The matrix $\langle L\rangle-L_*$ will then be positive semi-definite if and only if \refe{elemb}. holds and
\beq \label{elema}
a_*\geq a_L,
\eeq
\beq\label{elemcirc}
\dis \left({c_L\over a_L}-{c_*\over a_*}\right)^2 \dis \leq\left({1\over a_L}-{1\over a_*}\right)\left(d_L+{c^2_L\over a_L}-a_*-{c_*^2\over a_*}\right)
\eeq
By multiplying the last inequality by $a_*a_L$ and expanding (and using the fact that $a_*a_L>0$) we get \refe{simpcirc}.. Also the inequality \refe{elema}. is
superfluous as it is implied by \refe{simpcirc}. and the inequality $d_L\geq a_L$.
\section{Bounds on magneto-transport: 2d versus 3d}\label{s.cex}
From \refe{bou.abc}. and the non-negativity of the 4th (or 5th) diagonal element of $\langle L\rangle-L_*$ 
we deduce the additional (superfluous) bound
\beq\label{newbd}
2\,|c_*|\leq a_*+{c_*^2\over a_*}\leq\left\langle a+{c^2\over a}\right\rangle\leq\overline{a}+{\overline{c}^2\over\underline{a}},
\eeq
where to obtain the first inequality we have used the fact that $x+1/x\geq 2$ for all $x>0$. Thus, we have obtained an upper bound 
on $|c_*|$, but one that involves not only $\overline{c}$ but also $\overline{a}$ and $\underline{a}$. This is in contrast to the 
case for a columnar conductivity $\si(y)$ which is independent of the $y_3$-variable, where using the positivity of the 
determinant of the electric field $DU(y)=DU(y_1,y_2)$ established by Alessandrini and Nesi \cite{AlNe1,AlNe2}, we proved in \cite{BM3} 
that $c_*$ satisfies the same bound $\overline{c}$ as the local coefficient $c(y)$ (i.e., $|c_*|\leq\overline{c}$). 

Let us now relax the assumption that the effective tensor takes the form \refe{isocon*}.. 
The composite is said to be partially isotropic if the antisymmetric part of $\si_*$ satisfies
\beq\label{isop.sig*}
(\si_*)^{A}=c_*\,J,\quad\mbox{where}\quad
J:=\begin{pmatrix}
0 & -1 & 0
\\
1 & 0 & 0
\\
0 & 0 & 0
\end{pmatrix}.
\eeq
Given a partially isotropic composite we can always subdivide it into square columns with edges parallel to the $y_3$-axis and with side length
much larger than the existing microstructure, and then rotate each square column  about
its center axis by either $0^\circ$, $90^\circ$, $180^\circ$ or $270^\circ$ with equal probability in an 
uncorrellated way. The resulting polycrystal is invariant under rotations of $90^\circ$ about the $y_3$-axis and thus
will have an effective tensor of the form \refe{isocon*}.
and by a lemma of Stroud and Bergman \cite{StBer} will have the same constant $c_*$ as the original partially isotropic composite.

The question naturally arises as to whether for partially isotropic composites
$|c_*|$ can be bounded solely in terms of $\overline{c}$, like in the case of a columnar conductivity $\si(y)$?
The answer is no, it cannot. Indeed, we have the following result:

\begin{Thm}\label{thm.bou}
Consider a periodic conductivity $\si(y)$ given by \refe{iso.sig}. which satisfies the bounds \refe{bou.abc}.. Assume that the composite is partially isotropic in the sense of~\refe{isop.sig*}.. Then, the effective coefficient $c_*$ satisfies 
\beq\label{bouc*}
|c_*|\leq \frac{\overline{c}}{\underline{a}}\left(\si_* e_1\cdot e_1\right)^{1/2}\left(\si_* e_2\cdot e_2\right)^{1/2}.
\eeq
On the other hand, given any arbitrarily large constant $\ka>0$ there exist transversely isotropic conductivities $\si_{\th,\ka}(y)$ depending on a 
parameter $\th>0$, with $c(y)\in\{0,1\}$ a.e. $y\in\RR^3$, such that as $\th\to 0$ the effective conductivity $\si_{\th,\ka}^*$ is partially isotropic and satisfies
\beq\label{cthka*}
\lim_{\th\to 0}\,(\si_{\th,\ka}^*)^{A}=-\ka\,J\quad\mbox{or}\quad\lim_{\th\to 0}\,(\si_{\th,\ka}^*)^{A}=\ka\,J.
\eeq
\end{Thm}
\begin{Rem}\label{rem.isot}
In the case when $\si_*$ is transversely isotropic, taking the form \refe{isocon*}., the bound \refe{bouc*}. reduces to
\beq\label{newbdiso} 
|c_*|\leq\frac{a_*}{\underline{a}}\,\overline{c},
\eeq
and using \refe{newbd}. we obtain that
\beq |c_*|\leq \frac{\overline{c}}{\underline{a}}\left(\left\langle a+{c^2\over a}\right\rangle-{c_*^2\over a_*}\right) 
\leq \frac{\overline{c}}{\underline{a}}\left(\dis \overline{a}+{\overline{c}^2\over\underline{a}}\right).
\eeq
In contrast to the bounds \refe{newbd}. and \refe{simpcirc}. this new bound shows that $c_*$ necessarily goes to zero as $\overline{c}$ goes to zero if $\overline{a}$ and $\underline{a}$ are held fixed. Also 
if we add the antisymmetric matrix $c_0J$ to $\si(y)$ then the effective tensor will change to $\si_*+c_0J$, implying from \refe{newbdiso}. that the inequality
\beq |c_*+c_0|\leq \frac{a_*}{\underline{a}}\,\max\big(|c_++c_0|,|c_-+c_0|\big)
\eeq
holds for all constants $c_0$, where
\beq c_+:=\sup_{y\in Y}c(y),\quad c_-:=\inf_{y\in Y}c(y).
\eeq
Taking the optimum value $c_0=-(c_++c_-)/2$ gives the bounds
\beq\label{bestbd}
\left|2\,c_*-c_+-c_-\right|\leq \frac{a_*}{\underline{a}}\,(c_+-c_-).
\eeq 
\end{Rem}
\begin{Rem}\label{rem.npb}
Theorem~\ref{thm.bou} proves that contrary to the columnar case of \cite{BM3} we cannot expect to bound the effective coefficient $c_*$ only in terms of the bound $\overline{c}$ of the local coefficient $c(y)$. Actually, \refe{cthka*}. shows that arbitrarily large (positive or negative) effective coefficients $c_*$ can be derived although the local coefficient $c(y)$ only takes values in $\{0,1\}$. Here, the contrast of the symmetric part of the conductivity plays a crucial role as suggested in the bound  \refe{bouc*}.. This is strongly linked to the fact that the $(2\times 2)$ determinant
\beq\label{D12}
\De_{1,2}(DU):=\partial_1 u_1\,\partial_2 u_2-\partial_1 u_2\,\partial_2 u_1,\quad\mbox{for }U=(u_1,u_2,u_3),
\eeq
does not always have a constant sign throughout the material (see the proof of Theorem~\ref{thm.bou} below) contrary to the columnar case.
\end{Rem}
\noindent
{\bf Proof of Theorem~\ref{thm.bou}.}
\par\ms\noindent
{\it Bound for $c_*$:} The div-curl lemma of Murat-Tartar (see \cite{Tar1,Mur1,Mur2}) and the formula \refe{sig*}. for $\si_*$ yield
\beq\label{DUTsiDU}
\big\langle(DU)^T\si DU\big\rangle=\big\langle(DU)^T\big\rangle\,\big\langle\si DU\big\rangle
=\big\langle DU\rangle^T\si_*\big\langle DU\big\rangle=\si_*.
\eeq
Hence passing to the antisymmetric part it follows that
\beq
(\si_*)^{A}=\big\langle(DU)^T\si^{A} DU\big\rangle=\big\langle c\,(DU)^T JDU\big\rangle
=\left\langle c\left(\begin{smallmatrix}
0 & -\De_{1,2}(DU) & -\De_{1,3}(DU) \\�\De_{1,2}(DU) & 0 & -\De_{2,3}(DU) \\ \De_{1,3}(DU) & \De_{2,3}(DU) & 0
\end{smallmatrix}\right)\right\rangle,
\eeq
where $\De_{i,j}(DU):=\partial_1 u_i\,\partial_2 u_j-\partial_1 u_j\,\partial_2 u_i$. Therefore, since $\si_*$ is partially isotropic, we obtain the following formula for the effective coefficient $c_*$,
\beq\label{c*1}
c_*=\big\langle c\,\De_{1,2}(DU)\big\rangle.
\eeq
Using the Cauchy-Schwartz inequality we have
\beq\label{bo1}
\ba{ll}
|c_*| & \leq \overline{c}\left\langle 
\begin{pmatrix} 
|\partial_1 u_1|
\\
|\partial_2 u_1| 
\end{pmatrix}
\cdot
\begin{pmatrix} 
|\partial_2 u_2|
\\
|\partial_1 u_2| 
\end{pmatrix}
\right\rangle
\\ \ecart
& \leq \overline{c}\,\big\langle|\partial_1 u_1|^2+|\partial_2 u_1|^2\big\rangle^{1/2}\big\langle|\partial_1 u_2|^2+|\partial_2 u_2|^2\big\rangle^{1/2}.
\ea
\eeq 
On the other hand  \refe{DUTsiDU}. also implies  for $i=1,2$,
\beq\label{bo2}
\si_* e_i\cdot e_i= \big\langle\si\nabla u_i\cdot\nabla u_i\big\rangle\geq
\underline{a}\,\big\langle|\partial_1 u_i|^2+|\partial_2 u_i|^2\big\rangle.
\eeq
Combining \refe{bo1}. and \refe{bo2}. gives the desired bound \refe{bouc*}..
\par\ms\noindent
{\it Derivation of arbitrarily large coefficients $c_*$:} Let $\th$, $\ka$ be two positive numbers, let $\xi^1$, $\xi^2_\th$ be the vectors defined by
\beq\label{xii}
\xi^1_\th:=\left(0,{\th\over\sqrt{1+\th^2}},{1\over\sqrt{1+\th^2}}\right),\quad\xi^2:=\left(0,{1\over\sqrt{2}},{1\over\sqrt{2}}\right),
\eeq
and let $\si^1_{\th,\ka}$, $\si^2$, $\si^3$ be the (transversely isotropic) phases defined by
\beq\label{sigi}
\si^1_{\th,\ka}:=\begin{pmatrix} \ka\,\th^{-2} & 0 & 0 \\ 0 & \ka\,\th^{-2} & 0 \\ 0 & 0 & 1 \end{pmatrix},\quad\si^2:=I,\quad\si^3:=2\,I+J.
\eeq
Consider the rank-two laminate mixing in the direction $\xi^1_\th$ the phase $\si^1_{\th,\ka}$, with volume fraction $1-\th$, and the rank-one laminate, with volume fraction $\th$, composed of the mixture in the direction $\xi^2$ of the phases $\si^2$ and $\si^3$ with volume fraction ${1\over 2}$. The two-scale conductivity $\si_{\th,\ka}$ is defined by
\beq\label{sithka}
\si_{\th,\ka}(y,z):=\chi_\th(\xi^1_\th\cdot y)\,\si^1_{\th,\ka}+\big(1-\chi_\th(\xi^1_\th\cdot y)\big)
\Big(\chi(\xi^2\cdot z)\,\si^2+\left(1-\chi(\xi^2\cdot z)\right)\si^3\Big),
\eeq
where $y={x\over\ep}$, $z={x\over\ep^2}$ are the ordered fast variables, $\chi_\th$ is the $1$-periodic function which agrees with the characteristic function of $[0,1-\th]$ in $[0,1]$, and $\chi$ is the $1$-periodic function which agrees with the characteristic function of $[0,{1\over 2}]$ in $[0,1]$. By \cite{Mil86,Bri} the local electric field $E_{\th,\ka}$ associated with the conductivity $\si_{th,\ka}$ has the same laminate structure as \refe{sithka}., and thus can be written as
\beq\label{Ethka}
E_{\th,\ka}(y,z):=\chi_\th(\xi^1_\th\cdot y)\,E^1_{\th,\ka}+\left(1-\chi_\th(\xi^1_\th\cdot y)\right)
\left(\chi(\xi^2\cdot z)\,E^2_{\th,\ka}+\left(1-\chi(\xi^2\cdot z)\right)E^3_{\th,\ka}\right).
\eeq
The constant matrices $E^1_{\th,\ka}$, $E^2_{\th,\ka}$, $E^3_{\th,\ka}$ are the solutions of the linear system
\beq\label{Eithka}
\left\{\ba{ll}
(1-\th)\,E^1_{\th,\ka}+{\th\over 2}\,(E^2_{\th,\ka}+E^3_{\th,\ka})=I & \mbox{average-value}
\\ \ecart
E^2_{\th,\ka}-E^3_{\th,\ka}=\xi^2\otimes\eta_2 & \mbox{jump of the curl at the scale }\ep^2
\\
E^1_{\th,\ka}-{1\over 2}\,(E^2_{\th,\ka}+E^3_{\th,\ka})=\xi^1_\th\otimes\eta_1 & \mbox{jump of the curl at the scale }\ep
\\ \ecart
(\si^2 E^2_{\th,\ka}-\si^3 E^3_{\th,\ka})^T\xi^2=0 & \mbox{jump of the div at the scale }\ep^2
\\
 \big[\si^1_{\th,\ka} E^1_{\th,\ka}-{1\over 2}\,(\si^2 E^2_{\th,\ka}+\si^3 E^3_{\th,\ka})\big]^T\xi^1_\th=0
& \mbox{jump of the curl at the scale }\ep.
\ea\right.
\eeq
We refer to \cite{Bri} for more details. Similarly to \refe{DUTsiDU}. and taking into account the two-scale structure \refe{sithka}. the effective conductivity $\si^*_{\th,\ka}$ is given by
\beq\label{si*thka}
\si^*_{\th,\ka}=(1-\th)\,(E^1_{\th,\ka})^T\si^1_{\th,\ka} E^1_{\th,\ka}+{\th\over 2}
\left[(E^2_{\th,\ka})^T\si^2 E^2_{\th,\ka}+(E^3_{\th,\ka})^T\si^3 E^3_{\th,\ka}\right].
\eeq
Taking into account the values \refe{sigi}. of the matrix conductivities we deduce that
\beq
(\si^*_{\th,\ka})^{A}={\th\over 2}\,(E^3_{\th,\ka})^T J E^3_{\th,\ka}.
\eeq
Using Maple to compute explicitly the solutions $E^1_{\th,\ka}$, $E^2_{\th,\ka}$, $E^3_{\th,\ka}$ of the linear system \refe{Eithka}., we get the following asymptotics as $\th\to 0$,
\beq\label{negc*}
(\si^*_{\th,\ka})^{A}={\th\over 2}\,\De_{1,2}\left(E^3_{\th,\ka}\right)J+O(\th)=-{\ka\over 17}\,J+O(\th)
\eeq
Therefore, $\si^*_{\th,\ka}$ is asymptotically partially isotropic, and the effective coefficient
\beq\label{nc*thka}
c^*_{\th,\ka}:={\th\over 2}\,\De_{1,2}\left(E^3_{\th,\ka}\right)=-{\ka\over 17}+O(\th),
\eeq
is both negative and arbitrarily large when $\ka$ is arbitrarily large. Moreover, the $(2\times 2)$ determinant $\De_{1,2}$ of the electric field satisfies
\beq\label{nD12}
\De_{1,2}\left(E^3_{\th,\ka}\right)=-\De_{1,2}\left(E^2_{\th,\ka}\right)+O(1)=-{2\,\ka\over 17\,\th}+O(1),
\eeq
and thus for large $\ka$ has not the same sign throughout the material, contrary to the columnar case.
\par
On the other hand, if we replace in \refe{sigi}. the matrix $\si^3$ by
\beq\label{si3}
\si^3:=\begin{pmatrix} 2 & -1 & 0 \\ 1 & 2 & 0 \\ 0 & 0 & {1\over 2} \end{pmatrix},
\eeq
then the previous procedure leads us to the asymptotics
\beq\label{posc*}
(\si^*_{\th,\ka})^{A}={\th\over 2}\,\De_{1,2}\left(E^3_{\th,\ka}\right)J+O(\th)={\ka\over 13}\,J+O(\th).
\eeq
Hence, the effective conductivity $\si^*_{\th,\ka}$ is still asymptotically partially isotropic, and the effective coefficient
\beq\label{pc*thka}
c^*_{\th,\ka}={\th\over 2}\,\De_{1,2}\left(E^3_{\th,\ka}\right)={\ka\over 13}+O(\th),
\eeq
is arbitrarily large but positive. As before, the minor $\De_{1,2}$ of the electric field satisfies
\beq\label{pD12}
\De_{1,2}\left(E^3_{\th,\ka}\right)=-\De_{1,2}\left(E^2_{\th,\ka}\right)+O(1)={2\,\ka\over 13\,\th}+O(1),
\eeq
and does not have the same sign throughout the material when $\ka$ is large.
\par
\section{Hashin-Shtrikman type bounds under geometric isotropy}\label{s.HSbou}
\subsection{$Y$-tensors, $\Ga$-operator, and geometric isotropy}\label{s.Yisogeo}
For given $\alpha,\beta>0$, consider a periodic two-phase composite with local conductivity
\beq\label{si(y)}
\si(y)=\chi_1(y)\,\si_1+\chi_2(y)\,\si_2\in{\cal M}_\sharp(\alpha,\beta;Y),
\eeq
where $\chi_i$ is the characteristic function of the phase $i$ with volume fraction $f_i$, for $i=1,2$. Denote by $\si^*$ its effective conductivity. Following \cite{GiMi} (see also \cite{MilBook}, Chapter~19) there exists an effective tensor $Y_*$ associated with the conductivity $\si$, defined by ($e\in L^2_\sharp(Y)^3$ is the electric field and  $j\in L^2_\sharp(Y)^3$ is the current field)
\beq\label{PY*}
P(j)=-\,Y_*\,P(e),\quad\mbox{where}\quad j=\si\,e\quad\mbox{and}\quad P(g):=\big\langle\chi_1\left(g-\langle g\rangle\right)\big\rangle.
\eeq
In some sense $P$ is the projection on phase~$1$ of the fluctuating component of the field. Also we have for the adjoint problem
\beq\label{PY*T}
P(j')=-\,Y_*^TP(e'),\quad\mbox{where}\quad j'=\si^Te'.
\eeq
Recall the relation \refe{LEF}. and the definition \refe{L(y)}. of the tensor $L(y)$ which enters it. A similar
computation based on the formulas \refe{PY*}. and \refe{PY*T}. lead us to an effective tensor $\Y_*$ associated with the tensor $L(y)$ of \refe{L(y)}., and defined by
\beq\label{PcY*}
P\begin{pmatrix} e_S \\ j_A \end{pmatrix}=\begin{pmatrix} P(e_S) \\ P(j_A)\end{pmatrix}
=-\,\Y_*\begin{pmatrix} P(j_S) \\ P(e_A) \end{pmatrix}=-\,\Y_*\,P\begin{pmatrix} j_S \\ e_A \end{pmatrix},
\eeq
where similarly to \refe{L*}.,
\beq\label{cY*1}
\Y_*=\begin{pmatrix} (Y_*^S)^{-1} & (Y_*^S)^{-1}Y_*^A \\ \ecart -Y_*^A(Y_*^S)^{-1} & Y_*^S-Y_*^A(Y_*^S)^{-1}Y_*^A \end{pmatrix}.
\eeq
\par
Now, we will derive a Hashin-Shtrikman type variational inequality associated with the variational principle \refe{varpri}.. To this end, let us consider for a given reference tensor $L_0$, the nonlocal operator $\Ga$ defined for periodic vector-valued functions $A,B\in L^2_\sharp(Y)^6$, by
\beq\label{Ga}
B=\Ga A\qquad\mbox{if}\qquad\Ga_1 B=B\;\;\mbox{and}\;\;\Ga_1(A-L_0B)=A-L_0B,
\eeq
where $\Ga_1$ represents the projection on the space of fields which satisfy the same differential constraints as $E\in L^2_\sharp(Y)^6$ in \refe{LEF}.. Since $E$ is composed by a divergence free field $j_S\in L^2_\sharp(Y)^3$, and a curl free field $e_A\in L^2_\sharp(Y)^3$, the operator $\Ga_1$ in Fourier space is given by
\beq\label{Ga1xi}
\Ga_1(k)=\Ga_1(\xi)=\begin{pmatrix} I-\xi\otimes\xi & 0 \\ 0 & \xi\otimes\xi \end{pmatrix},\quad\mbox{where}\quad\xi:={k\over|k|},
\quad\mbox{ for }k\in\ZZ^3\setminus\{0\}.
\eeq
Under the conditions $L_i>L_0\geq 0$, for $i=1,2$, the Hashin-Shtrikman type variational inequality associated with the variational principle \refe{varpri}. is given by the formula~(13.30) of \cite{MilBook}, which reads as
\beq\label{HS1}
\left(L_*-L_0\right)^{-1}\langle F\rangle:\langle F\rangle\leq\big\langle\left[\Ga+(L_*-L_0)^{-1}\right]F:F\big\rangle,
\quad\mbox{for any }F\in L^2_\sharp(Y)^6.
\eeq
Following the computations of \cite{MilBook} (Section~23.6) this inequality implies the bound
\beq\label{HS2}
\Y_*+L_0\geq\left[{1\over f_1f_2}\sum_{k\in\ZZ^3\setminus\{0\}}\big|{\hat{\chi}_1(k)}\big|^2\,\Ga(k)\right]^{-1},
\eeq
which also holds for the enlarged inequalities $L_i\geq L_0\geq 0$, for $i=1,2$.
Note that by virtue of the Plancherel equality the Fourier coefficients $\hat{\chi}_1(k)$ of the characteristic function $\chi_1$ satisfy the equality
\beq
{1\over f_1f_2}\sum_{k\in\ZZ^3\setminus\{0\}}\big|\hat{\chi}_1(k)\big|^2={1\over f_1f_2}\,\big\langle(\chi_1-f_1)^2\big\rangle=1.
\eeq
So, the series in \refe{HS2}. can be regarded of an average of the operator $\Ga$.
\par
Finally, consider the case of a two-phase random composite. According to \cite{Wil} (see also \cite{MilBook}, Section~15.6) the composite is said to have a geometric isotropy if all correlation functions associated with the geometry represented by the characteristic function $\chi_1$ are invariant by rotation (or reflection). Then, under geometric isotropy the average series  of \refe{HS2}. reduces to an average of $\Ga$ over all directions of the unit sphere $S^2$. Therefore, we get the bound (see \cite{MilBook}, Section~23.6) 
\beq\label{HS}
\Y_*+L_0\geq{1\over\langle\Ga\rangle_{S^2}}.
\eeq
\subsection{Hashin-Shtrikman type bounds}
\begin{figure}[!t]
\centering
\includegraphics[scale=.5]{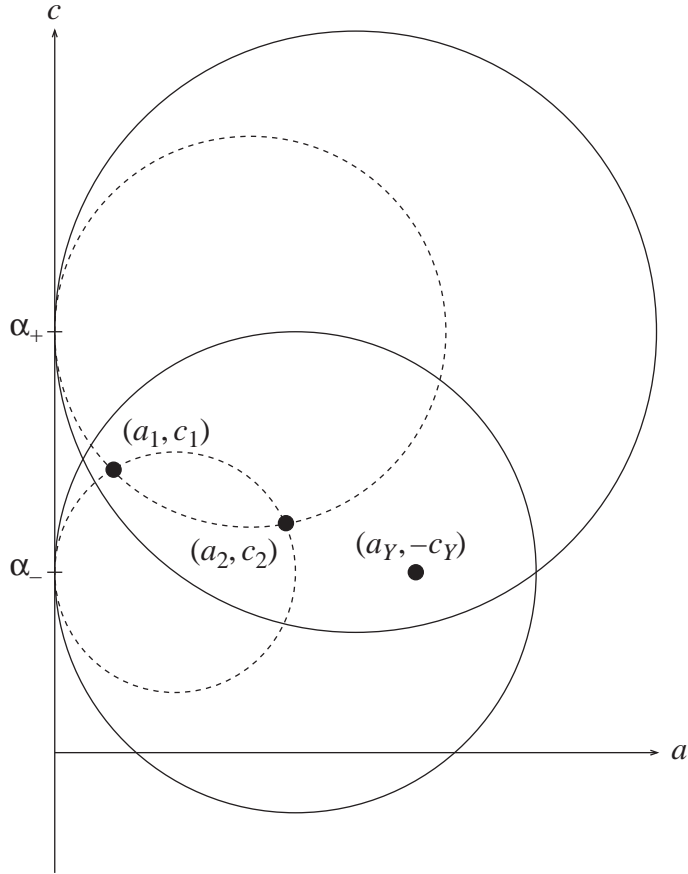}
\caption{The circles $\bigcirc^{\mbox{\tiny HS}}_\pm$ surrounding $(a_Y,-c_Y)$ and the (dashed) circles $\bigcirc_\pm$ passing through $(a_i,c_i)$, for $i=1,2$, assuming $\al_-\leq\al_+$.}
\label{fig1}
\end{figure}

Consider a periodic two-phase composite with non-symmetric positive definite conductivities
\beq\label{sig12}
\si_i:=\begin{pmatrix} a_i & -c_i & 0 \\ c_i & a_i & 0 \\ 0 & 0 & b_i \end{pmatrix},\quad\mbox{with }b_1\geq b_2,
\eeq
of respective volume fractions $f_i$, for $i=1,2$. The effective conductivity $\si_*$ of the composite is assumed to be transversely isotropic, i.e.
\beq\label{sigma*}
\si_*:=\begin{pmatrix} a_* & -c_* & 0 \\ c_* & a_* & 0 \\ 0 & 0 & b_* \end{pmatrix}.
\eeq
Let $g:(0,\infty)\to\RR$ be the function defined by
\beq\label{g}
g(r):={1\over 2}\int_0^\pi{\cos^2\th\,\sin\th\over \cos^2\th+r^{-1}\sin^2\th}\,d\th\in(0,1),\quad\mbox{for }r>0.
\eeq
Consider the coefficients $\al_\pm$, $t_1^\pm$, $s_1^\pm$, $a_Y$, $b_Y$, $c_Y$ defined by
\beq\label{alpm}
\al_\pm:={a_1 c_2-a_2 c_1\pm\sqrt{a_1a_2\,\big((a_1-a_2)^2+(c_1-c_2)^2\big)}\over a_1-a_2},
\eeq
\beq\label{sipm}
t_1^\pm:={a_1\over a_1^2+(c_1-\al_\pm)^2},\quad
s_1^\pm:={2\,t_1^\pm\over 1+g(b_1t_1^\pm)}-t_1^\pm,
\eeq
\beq\label{acY}
a_Y+i\,c_Y:=-f_2(a_1+i\,c_1)-f_1(a_2+i\,c_2)
+{f_1f_2\left(a_1+i\,c_1-a_2-i\,c_2\right)^2\over\big(f_1(a_1+i\,c_1)+f_2(a_2+i\,c_2)-a_*-i\,c_*\big)},
\eeq
\beq\label{bY}
b_Y:=-f_2 b_1-f_1 b_2+{f_1f_2\left(b_1-b_2\right)^2\over f_1 b_1+f_2 b_2-b_*}.
\eeq
Then, we have the following result:
\begin{Thm}\label{thm.HSbou}
Assume that the composite is geometrically isotropic.
Then, in view of definitions \refe{alpm}.-\refe{bY}. the coefficients $a_*$, $c_*$ of the effective conductivity $\si_*$ \refe{sigma*}. satisfy the Hashin-Shtrikman type bounds
\beq\label{bouac*}
a_Y^2+(c_Y+\al_\pm)^2-{a_Y\over s_1^\pm}\leq 0,
\eeq
while the points $(a_1,c_1)$, $(a_2,c_2)$ solve the equations
\beq\label{bouac}
a^2+(c-\al_\pm)^2-{a\over t_1^\pm}=0.
\eeq
Moreover, the coefficient $b_*$ satisfies the bounds
\beq\label{boub*}
{1\over b_Y}+{1\over b_1}\geq {1\over b_1\big(1-g(b_1t_1^\pm)\big)},\quad b_Y\geq 0.
\eeq
\end{Thm}
\begin{Rem}\label{rem.HSbou1}
In the $a_Y$-$c_Y$ plane the bounds \refe{bouac*}. correspond to the intersection of two disks parametrized by $\al_\pm$, which are tangent to the $c_Y$-axis. Due to definition \refe{acY}. these bounds remain the same if we replace $c_1$, $c_2$, $-c_Y$ by $c_1+c_0$, $c_2+c_0$, $-c_Y+c_0$. This reflects the fact if we add a  antisymmetric matrix to the local conductivity $\si$, then the same antisymmetric matrix is added to~$\si_*$. Also note that if $c_1=c_2=c_*=c$, then $c_Y=-c$.
\end{Rem}
\begin{Rem}\label{rem.HSbou2}
With the change $c_Y$ to $-c_Y$, the circle $\bigcirc^{\mbox{\tiny HS}}_\pm$ satisfying the equality in \refe{bouac*}. is the same as the circle $\bigcirc_\pm$ of equation \refe{bouac}. passing through the points $(a_i,c_i)$, $i=1,2$, when $s_1^\pm=t_1^\pm$. Moreover, since $g(r)\in(0,1)$ for $r>0$, we have $0<s_1^\pm<t_1^\pm$ in \refe{sipm}.. This implies that the radius $(2s_1^\pm)^{-2}$ of $\bigcirc^{\mbox{\tiny HS}}_\pm$ is larger than the radius $(2t_1^\pm)^{-2}$ of $\bigcirc_\pm$. The circles $\bigcirc^{\mbox{\tiny HS}}_\pm$ and $\bigcirc_\pm$ are also tangent at the same point of the $c\,$-axis. The geometrical picture is given by Figure~\ref{fig1} in the $a$-$c$ plane.
\end{Rem}

\begin{Rem}\label{rem.HSbou3}
The inequalities \refe{bouac*}., \refe{boub*}. do not allow us to show that $c_*$ tends to zero when $c_1$ and $c_2$ approach zero, while keeping $a_1$, $a_2$, $b_1$, and $b_2$ fixed. 
To this end, we only have the bound \refe{bestbd}. which reads as
\beq\label{bestbd2}
|2\,c_*-c_1-c_2|\leq \frac{a_*}{\min\left(a_1,a_2\right)}\,|c_1-c_2|.
\eeq 
\end{Rem}
\noindent
{\bf Proof of Theorem~\ref{thm.HSbou}.}
The proof is divided in four steps. In the first step we determine a suitable reference tensor $L_0$. In the second step we compute the tensor $\Ga(\xi)$ involved in the $Y$-tensor approach. In the third step we compute the average $\langle\Ga\rangle_{S^2}$. 
The fourth step is devoted to the derivation of the bounds.
\par\ms\noindent
{\it First step~}: Determination of $L_0$.
\par\noindent
Similarly to \refe{L(y)}., let $L_i$, for $i=1,2,*$, be the symmetric tensor defined by
\beq\label{Li}
L_i:=\begin{pmatrix}
(\si_i^S)^{-1} & -\,(\si_i^S)^{-1}\,\si_i^A
\\ \ecart
\si_i^A\,(\si_i^S)^{-1} & \si_i^S-\si_i^A\,(\si_i^S)^{-1}\,\si_i^A
\end{pmatrix}
=\begin{pmatrix}
\dis {1\over a_i} & 0 & 0 & 0 & \dis {c_i\over a_i} & 0
\\
0 & \dis {1\over a_i} & 0 & \dis -\,{c_i\over a_i} & 0 & 0
\\
0 & 0 & \dis {1\over b_i} & 0 & 0 & 0
\\
0 & \dis -\,{c_i\over a_i} & 0 & \dis a_i+{c_i^2\over a_i} & 0 & 0
\\
\dis {c_i\over a_i} & 0 & 0 & 0 & \dis a_i+{c_i^2\over a_i} & 0
\\
0 & 0 & 0 & 0 & 0 & b_i
\end{pmatrix}.
\eeq
Now, let $L_0$ be the symmetric tensor defined by
\beq\label{L0}
L_0:=\begin{pmatrix}
C_1 & C_2
\\ \ecart
C_2^T & C_3
\end{pmatrix}
=\begin{pmatrix}
\dis t_1 & 0 & 0 & 0 & t_2 & 0
\\
0 & t_1 & 0 & -\,t_2 & 0 & 0
\\
0 & 0 & t_4 & 0 & 0 & 0
\\
0 & -\,t_2 & 0 & t_3 & 0 & 0
\\
t_2 & 0 & 0 & 0 & t_3 & 0
\\
0 & 0 & 0 & 0 & 0 & t_5
\end{pmatrix},
\quad\mbox{where }C_j\in\RR^{3\times 3},\ C_2^T=-\,C_2.
\eeq
The condition $L_0\geq 0$ is equivalent to
\beq\label{ti}
t_1\geq 0,\quad t_4\geq 0,\quad t_3\geq 0, \quad t_5\geq 0,\quad t_1 t_3\geq t_2^2\,.
\eeq
We also need $L_i\geq L_0$, which is equivalent to
\beq\label{inet1t2}
{1\over a_i}\geq t_1,\quad\det\begin{pmatrix}\dis {1\over a_i}-t_1 & \dis {c_i\over a_i}-t_2
\\ \ecart
\dis {c_i\over a_i}-t_2 & \dis a_i+{c_i^2\over a_i}-t_3\end{pmatrix}\geq 0,\quad\mbox{for }i=1,2,
\eeq
\beq\label{inet4t5}
{1\over b_i}\geq t_4,\quad b_i\geq t_5,\quad\mbox{for }i=1,2.
\eeq
From now on, assume that $b_1\geq b_2$, and set
\beq\label{t4t5}
t_4:={1\over b_1},\quad t_5=b_2,
\eeq
in order to make the  inequalities \refe{inet4t5}. as sharp as possible.
\par
On the other hand, the inequalities \refe{inet1t2}. show that the points $(a_1,c_1)$ and $(a_2,c_2)$ belong to the disk in the $a$-$c$ plane,
\beq\label{acti}
t_1\,(a^2+c^2)+(t_2^2-t_1t_3-1)\,a-2t_2\,c+t_3\leq 0,
\eeq
which lies in the half-plane $a\geq 0$. To make these bounds as tight as possible, we consider the two circles which are tangent to the  $c$-axis, and which pass through the points $(a_1,c_1)$ and $(a_2,c_2)$. This requires
\beq\label{t3=}
t_3={t_2^2\over t_1},
\eeq
and the two circle equations
\beq\label{t1t2}
t_1^2\,(a_i^2+c_i^2)-t_1\,a_i-2t_1t_2\,c_i+t_2^2=0,\quad\mbox{for }i=1,2,
\eeq
which can be written as
\beq
\left\{\ba{l}
t_1^2\,a_2(a_1^2+c_1^2)-t_1\,a_1 a_2-2t_1t_2\,a_2 c_1+t_2^2\,a_2=0
\\ \ecart
t_1^2\,a_1(a_2^2+c_2^2)-t_1\,a_1 a_2-2t_1t_2\,a_1 c_2+t_2^2\,a_1=0.
\ea\right.
\eeq
Subtracting and dividing by $t_1^2$, we get that $\al:=t_2/t_1$ solves
\beq\label{equal}
(a_1-a_2)\,\al^2-2\,(a_1c_2-a_2c_1)\,\al+a_1(a_2^2+c_2^2)-a_2(a_1^2+c_1^2)=0,
\eeq
the discriminant of which is
\beq
(a_1 c_2-a_2 c_1)^2+\left(a_1-a_2\right)\big(a_2(a_1^2+c_1^2)-a_1(a_2^2+c_2^2)\big)
=a_1a_2\,\big((a_1-a_2)^2+(c_1-c_2)^2\big)\geq 0.
\eeq
Hence, equation \refe{equal}. has two real solutions (one for each circle) $\al_\pm$ which are given by \refe{alpm}..
Moreover, putting $t_2=\al\,t_1$ in \refe{t1t2}. we obtain that
\beq\label{t1=}
t_1={a_i\over a_i^2+(c_i-\al)^2}\leq{1\over a_i},\quad\mbox{for }i=1,2,
\eeq
which implies that the ($2\times 2$) matrix in \refe{inet1t2}. is non-negative. Therefore, the choice of the coefficients $t_1,t_2,t_3$ given by
\beq\label{t1t2t3}
t_1={a_1\over a_1^2+(c_1-\al)^2},\quad t_2=\al\,t_1,\quad t_3=\al^2\,t_1,\quad\mbox{for }\al=\al_\pm,
\eeq
combined with \refe{t4t5}., implies the desired inequalities $L_i\geq L_0$, for $i=1,2$. Making this choice in \refe{acti}. the points $(a_1,c_1)$ and $(a_2,c_2)$ belong to the two circles of equation \refe{bouac}. which are  tangent to the line $a=0$.
\par\ms\noindent
{\it Second step~}: Computation of $\Ga(\xi)$.
\par\noindent
Let $\xi\in S^2$, $\xi=(\sin\th \cos\ph,\sin\th \sin\ph,\cos\th)$. By virtue of Section~\ref{s.Yisogeo} the tensor $\Ga(\xi)$ is defined from the tensor $L_0$ \refe{L0}., by
\beq\label{Ga0}
\begin{pmatrix} B_1 \\ B_2 \end{pmatrix}=\Ga(\xi)\begin{pmatrix} A_1 \\ A_2 \end{pmatrix},\quad\mbox{for }A_1,A_2,B_1,B_2\in\RR^{3\times 3},
\eeq
if and only if
\beq\label{Ga1}
\begin{pmatrix} I-\xi\otimes\xi & 0 \\ 0 & \xi\otimes\xi \end{pmatrix}\begin{pmatrix} B_1 \\ B_2 \end{pmatrix}=\begin{pmatrix} B_1 \\ B_2 \end{pmatrix},
\eeq
\beq\label{Ga2}
\mbox{and}\quad\begin{pmatrix} I-\xi\otimes\xi & 0 \\ 0 & \xi\otimes\xi \end{pmatrix}\left[\begin{pmatrix} A_1 \\ A_2 \end{pmatrix}
-\begin{pmatrix} C_1 & C_2 \\ C_2^T & C_3 \end{pmatrix} \begin{pmatrix} B_1 \\ B_2 \end{pmatrix}\right]=0.
\eeq
By \refe{Ga1}. we have $B_1^T\xi=0$, and $B_2=\xi\otimes\eta$ for some vector $\eta$. From \refe{Ga2}. it follows that 
\beq\label{Ga3}
A_2^T\xi-B_1^TC_2\xi-B_2^TC_3^T\xi=A_2^T\xi-B_1^TC_2\xi-(C_3\xi\cdot\xi)\,\eta=0,
\eeq
\beq\label{B1}
\ba{ll}
A_1-C_1B_1-C_2\,(\xi\otimes\eta) & =\xi\otimes\left(A_1^T\xi-B_1^TC_1^T\xi-B_2^TC_2^T\xi\right)
\\ \ecart
& =\xi\otimes\left(A_1^T\xi-B_1^TC_1^T\xi-(C_2\xi\cdot\xi)\,\eta\right)
\\ \ecart
& =\xi\otimes\left(A_1^T\xi-B_1^TC_1^T\xi\right)\quad\mbox{since }C_2^T=-\,C_2,
\\ \ecart
& =\xi\otimes k\quad\mbox{where }k:=A_1^T\xi-B_1^TC_1^T\xi.
\ea
\eeq
Noting that $C_1^{-1}C_2$ is antisymmetric, this implies that
\beq
0=B_1^T\xi=(C_1^{-1}A_1)^T\xi-(\eta\otimes\xi)\,(C_1^{-1}C_2)^T\xi-(k\otimes\xi)\,(C_1^{-1})^T\xi
=(C_1^{-1}A_1)^T\xi-(C_1^{-1}\xi\cdot\xi)\,k,
\eeq
so we have
\beq\label{k}
k={(C_1^{-1}A_1)^T\xi\over C_1^{-1}\xi\cdot\xi}.
\eeq
Moreover, replacing $B_1$ given by \refe{B1}. in \refe{Ga3}. and using that $(C_1^{-1})^TC_2$ is antisymmetric, we get that
\beq
\ba{rl}
0 \kern -.5em & =A_2^T\xi-(C_1^{-1}A_1)^TC_2\xi+(\eta\otimes\xi)\,(C_1^{-1}C_2)^TC_2\xi+(k\otimes\xi)\,(C_1^{-1})^TC_2\xi-(\eta\otimes\xi)\,C_3^T\xi
\\ \ecart
& = A_2^T\xi-(C_1^{-1}A_1)^TC_2\xi-\left[(C_3-C_2^TC_1^{-1}C_2)\,\xi\cdot\xi\right]\eta,
\ea
\eeq
hence
\beq\label{etaD}
\eta={A_2^T\xi-(C_2^TC_1^{-1}A_1)^T\xi\over D\xi\cdot\xi},\quad\mbox{where }D:=C_3-C_2^TC_1^{-1}C_2.
\eeq
Again using \refe{B1}. combined with \refe{k}. and \refe{etaD}. we deduce that
\beq
\ba{rl}
B_1= & \dis C_1^{-1}A_1-{C_1^{-1}C_2\,(\xi\otimes\xi)\over D\xi\cdot\xi}\left(A_2-C_2^TC_1^{-1}A_1\right)-{C_1^{-1}\,(\xi\otimes\xi)\over C_1^{-1}\xi\cdot\xi}\,C_1^{-1}A_1
\\ \ecart
B_2= & \dis {\xi\otimes\xi\over D\xi\cdot\xi}\,A_2-{\xi\otimes\xi\over D\xi\cdot\xi}\,C_2^TC_1^{-1}A_1.
\ea
\eeq
Hence, from definition \refe{Ga0}. it follows that
\beq\label{Gaxi}
\Ga(\xi)=\begin{pmatrix}
\dis C_1^{-1}+{C_1^{-1}C_2\,(\xi\otimes\xi)\,C_2^TC_1^{-1}\over D\xi\cdot\xi}-{C_1^{-1}\,(\xi\otimes\xi)\,C_1^{-1}\over C_1^{-1}\xi\cdot\xi}
&& \dis {C_1^{-1}C_2^T(\xi\otimes\xi)\over D\xi\cdot\xi}
\\ \ecart
\dis {(\xi\otimes\xi)\,C_2C_1^{-1}\over D\xi\cdot\xi} && \dis {\xi\otimes\xi\over D\xi\cdot\xi}
\end{pmatrix}
\eeq
which is a symmetric matrix since $C_1^T=C_1$ and $C_2^T=-\,C_2$.
\par\ms\noindent
{\it Third step~}: Computation of $\left(\langle\Ga\rangle_{S^2}\right)^{-1}$.
\par\noindent
Note that the computation of $\Ga(\xi)$ can be carried out if the matrix $D$ of \refe{etaD}.
\beq\label{D}
D=\begin{pmatrix} d_1 & 0 & 0 \\ 0 & d_1 & 0 \\ 0 & 0 & d_2 \end{pmatrix}=\begin{pmatrix} \dis t_3-{t_2^2\over t_1} & 0 & 0 \\ 0 & \dis t_3-{t_2^2\over t_1} & 0 \\ 0 & 0 & t_5 \end{pmatrix},
\eeq
is positive definite, i.e. $t_1t_3>t_2^2$ and $t_5>0$. Let us assume these conditions for the moment.
We shall be able pass to the limit as $d_1\to 0$ in the expression of $\left(\langle\Ga\rangle_{S^2}\right)^{-1}$.
Set
\beq\label{P}
P=\begin{pmatrix} p_1 & 0 & 0 \\ 0 & p_1 & 0 \\ 0 & 0 & p_2 \end{pmatrix}:=\left\langle{\xi\otimes\xi\over D\xi\cdot\xi}\right\rangle_{S^2},
\eeq
\beq\label{QR}
Q=\begin{pmatrix} q_1 & 0 & 0 \\ 0 & q_1 & 0 \\ 0 & 0 & q_2 \end{pmatrix}:=\left\langle{\xi\otimes\xi\over C_1^{-1}\xi\cdot\xi}\right\rangle_{S^2},\quad
R=\begin{pmatrix} r_1 & 0 & 0 \\ 0 & r_1 & 0 \\ 0 & 0 & r_2 \end{pmatrix}:=C_1^{-1}-C_1^{-1}QC_1^{-1}.
\eeq
By definition \refe{moyf(xi)}. we have
\beq\label{p1}
p_1=\left\langle{\xi_1^2\over d_1+(d_2-d_1)\,\xi_3^2}\right\rangle_{S^2}\;\mathop{\longrightarrow}_{d_1\to 0}\;{1\over 4}\int_0^\pi{\sin^3\th\over d_2\,\cos^2\th}\,d\th=\infty,
\eeq
\beq\label{p2}
p_2=\left\langle{\xi_3^2\over d_1+(d_2-d_1)\,\xi_3^2}\right\rangle_{S^2}\;\mathop{\longrightarrow}_{d_1\to 0}\;{1\over d_2}={1\over t_5}.
\eeq
Moreover, the matrix
\beq\label{C1Q}
C_1^{-{1\over 2}}\,Q\,C_1^{-{1\over 2}}=\begin{pmatrix} \dis {q_1\over t_1} & 0 & 0 \\ 0 & \dis {q_1\over t_1} & 0 \\ 0 & 0 & \dis {q_2\over t_4} \end{pmatrix}
=\left\langle{C_1^{-{1\over 2}}\,(\xi\otimes\xi)\,C_1^{-{1\over 2}}\over C_1^{-1}\xi\cdot\xi}\right\rangle_{S^2}
\eeq
has the property that its trace is $1$. This combined with definitions \refe{moyf(xi)}. and \refe{g}. yields
\beq\label{qi}
\left\{\ba{rl}
q_2= & \dis t_4\left\langle{\xi_3^2\over (t_4/t_1)(\xi_1^2+\xi_2^2)+\xi_3^2}\right\rangle_{S^2}=t_4\,g(t_1/t_4)
\\ \ecart
q_1= & \dis {t_1\over 2}\left(1-{q_2\over t_4}\right)={t_1\over 2}\big(1-g(t_1/t_4)\big),
\ea\right.
\eeq
which also implies that
\beq\label{ri}
\left\{\ba{rl}
r_1= & \dis {1\over t_1}-{q_1\over t_1^2}={1\over 2t_1}\big(1+g(t_1/t_4)\big)
\\ \ecart
r_2= & \dis {1\over t_4}-{q_2\over t_4^2}={1\over t_4}\big(1-g(t_1/t_4)\big).
\ea\right.
\eeq
On the other hand, by definition \refe{P}. we have
\beq\label{CiP}
C_1^{-1}C_2\,P\,C_2^TC_1^{-1}=\begin{pmatrix} \dis {t_2^2\over t_1^2}\,p_1 & 0 & 0 \\ 0 & \dis {t_2^2\over t_1^2}\,p_1 & 0 \\ 0 & 0 & 0 \end{pmatrix},\quad
C_1^{-1}C_2^TP=\begin{pmatrix} \dis 0 & \dis -\,{t_2\over t_1}\,p_1 & 0 \\ \dis {t_2\over t_1}\,p_1 & 0 & 0 \\ 0 & 0 & 0 \end{pmatrix}.
\eeq
Then, putting \refe{QR}. and \refe{CiP}. in the $S^2$-average of \refe{Gaxi}. we get that
\beq\label{moyGaxi}
\langle\Ga\rangle_{S^2}=\begin{pmatrix}
\dis r_1+{t_2^2\over t_1^2}\,p_1 & 0 & 0 & 0 & \dis -\,{t_2\over t_1}\,p_1 & 0
\\
0 & \dis r_1+{t_2^2\over t_1^2}\,p_1 & 0 & \dis {t_2\over t_1}\,p_1 & 0 & 0
\\
0 & 0 & r_2 & 0 & 0 & 0
\\
0 & \dis {t_2\over t_1}\,p_1 & 0 & p_1 & 0 & 0
\\
\dis -{t_2\over t_1}\,p_1 & 0 & 0 & 0 & p_1 & 0
\\
0 & 0 & 0 & 0 & 0 & p_2
\end{pmatrix},
\eeq
which gives
\beq
{1\over\langle\Ga\rangle_{S^2}}=\begin{pmatrix}
\dis {1\over r_1} & 0 & 0 & 0 & \dis {t_2\over t_1r_1} & 0
\\
0 & \dis {1\over r_1} & 0 & \dis -\,{t_2\over t_1r_1} & 0 & 0
\\
0 & 0 & \dis {1\over r_2} & 0 & 0 & 0
\\
0 & \dis -\,{t_2\over t_1r_1} & 0 & \dis {1\over p_1}+{t_2^2\over t_1^2\,r_1} & 0 & 0 
\\
\dis {t_2\over t_1r_1} & 0 & 0 & 0 & \dis {1\over p_1}+{t_2^2\over t_1^2\,r_1} & 0
\\
0 & 0 & 0 & 0 & 0 & \dis {1\over p_2}
\end{pmatrix}
\eeq
Therefore, passing to the limit as $d_1\to 0$, or equivalently $t_3\to t_2^2/t_1$, \refe{p1}. and \refe{p2}. imply that
\beq\label{moyGaxi1}
{1\over\langle\Ga\rangle_{S^2}}=\begin{pmatrix}
\dis {1\over r_1} & 0 & 0 & 0 & \dis {t_2\over t_1r_1} & 0
\\
0 & \dis {1\over r_1} & 0 & \dis -\,{t_2\over t_1r_1} & 0 & 0
\\
0 & 0 & \dis {1\over r_2} & 0 & 0 & 0
\\
0 & \dis -\,{t_2\over t_1r_1} & 0 & \dis {t_2^2\over t_1^2\,r_1} & 0 & 0 
\\
\dis {t_2\over t_1r_1} & 0 & 0 & 0 & \dis {t_2^2\over t_1^2\,r_1} & 0
\\
0 & 0 & 0 & 0 & 0 & \dis t_5
\end{pmatrix}.
\eeq
\par\ms\noindent
{\it Fourth step~}: Derivation of the bounds.
\par\noindent
On the one hand, the Appendix of \cite{GiMi} (see also formula (19.3) of \cite{MilBook}) yields the following formula for the $Y$-tensor defined by \refe{PY*}.
\beq\label{Y*si*}
Y_*=-f_2\,\si_1-f_1\,\si_2+f_1f_2\left(\si_1-\si_2\right)\left(f_1\,\si_1+f_2\,\si_2-\si_*\right)^{-1}\left(\si_1-\si_2\right).
\eeq
Note that, due to the transverse isotropy of $\si_i$, for $i=1,2,*$, we have
\beq\label{Y*}
Y_*=\begin{pmatrix} a_Y & -c_Y & 0 \\ c_Y & a_Y & 0 \\ 0 & 0 & b_Y \end{pmatrix}.
\eeq
This relation also separates into blocks, so that we obtain for the $_{33}$ entry of $Y_*$ the relation \refe{bY}..
Moreover, making the correspondence
\beq
\begin{pmatrix} a & -c \\ c & a \end{pmatrix}\;\longleftrightarrow\, a+i\,c,
\eeq
we deduce from the first $(2\times 2)$ block of \refe{Y*si*}. the relation \refe{acY}..
\par
On the other hand, by \refe{Y*}. the formula \refe{cY*1}. for $\Y^*$ reads as
\beq\label{cY*2}
\Y_*=\begin{pmatrix}
\dis {1\over a_Y} & 0 & 0 & 0 & \dis -\,{c_Y\over a_Y} & 0
\\
0 & \dis {1\over a_Y} & 0 & \dis {c_Y\over a_Y} & 0 & 0
\\
0 & 0 & \dis {1\over b_Y} & 0 & 0 & 0
\\
0 & \dis {c_Y\over a_Y} & 0 & \dis a_Y+{c_Y^2\over a_Y} & 0 & 0
\\
\dis -\,{c_Y\over a_Y} & 0 & 0 & 0 & \dis a_Y+{c_Y^2\over a_Y} & 0
\\
0 & 0 & 0 & 0 & 0 & b_Y
\end{pmatrix}.\eeq
Then, the bound \refe{HS}. applied with the formulas \refe{cY*2}. for $\Y_*$, \refe{L0}. for $L_0$ and \refe{moyGaxi1}. for $\left(\langle\Ga\rangle_{S^2}\right)^{-1}$, combined with \refe{t4t5}., \refe{qi}., and \refe{ri}., implies that
\beq\label{bYt1}
{1\over b_Y}+{1\over b_1}\geq{1\over r_2}={1\over b_1\big(1-\dis g(b_1t_1)\big)},\quad b_Y\geq 0,
\eeq
and
\beq\label{acYsi}
\det\begin{pmatrix} \dis {1\over a_Y}-s_1 & \dis {c_Y\over a_Y}+s_2
\\ \ecart
\dis {c_Y\over a_Y}+s_2 & \dis a_Y+{c_Y^2\over a_Y}-s_3\end{pmatrix}\geq 0,
\eeq
where
\beq\label{si}
\left\{\ba{rll}
s_1:= & \dis {1\over r_1}-t_1& \dis =t_1\left({2\over 1+\dis g(b_1t_1)}-1\right)\geq 0\quad(\mbox{since }0< g(b_1t_1)<1)
\\ \ecart
s_2:= & \dis {t_2\over t_1r_1}-t_2 & \dis =t_2\left({2\over 1+\dis g(b_1t_1)}-1\right)
\\ \ecart
s_3:= & \dis {t_2^2\over t_1^2r_1}-t_3.
\ea\right.
\eeq
Due to \refe{t3=}. we have $s_1s_3=s_2^2$. Therefore, similarly to \refe{inet1t2}. and \refe{acti}. the inequality \refe{acYsi}. can be written as
\beq\label{acYs1s2}
a_Y^2+\left(c_Y+{s_2\over s_1}\right)^2-{a_Y\over s_1}\leq 0.
\eeq
Finally, taking into account \refe{sipm}., \refe{t1t2t3}., \refe{si}. the inequalities \refe{acYs1s2}. and \refe{bYt1}. correspond respectively to the desired bounds \refe{bouac*}., \refe{boub*}.. Theorem~\ref{thm.HSbou} is proved. $\cqfd$
\par\bs\noindent
{\bf Acknowledgements.} GWM is grateful for support from the Mathematical Sciences Research Institute and from
National Science Foundation through grant DMS-0707978.

\end{document}